\documentclass{amsart} 
\usepackage{amssymb}
\usepackage{amsmath}
\usepackage{amsfonts}

\begin{document}
\newtheorem{theo}{Theorem}[section]
\newtheorem{prop}[theo]{Proposition}
\newtheorem{lemma}[theo]{Lemma}
\newtheorem{exam}[theo]{Example}
\newtheorem{coro}[theo]{Corollary}
\theoremstyle{definition}
\newtheorem{defi}[theo]{Definition}
\newtheorem{rem}[theo]{Remark}


\newcommand{\Bb}{{\bf B}}
\newcommand{\Nb}{{\bf N}}
\newcommand{\Qb}{{\bf Q}}
\newcommand{\Rb}{{\bf R}}
\newcommand{\Zb}{{\bf Z}}
\newcommand{\Ac}{{\mathcal A}}
\newcommand{\Bc}{{\mathcal B}}
\newcommand{\Cc}{{\mathcal C}}
\newcommand{\Dc}{{\mathcal D}}
\newcommand{\Fc}{{\mathcal F}}
\newcommand{\Ic}{{\mathcal I}}
\newcommand{\Jc}{{\mathcal J}}
\newcommand{\Lc}{{\mathcal L}}
\newcommand{\Oc}{{\mathcal O}}
\newcommand{\Pc}{{\mathcal P}}
\newcommand{\Sc}{{\mathcal S}}
\newcommand{\Tc}{{\mathcal T}}
\newcommand{\Uc}{{\mathcal U}}
\newcommand{\Vc}{{\mathcal V}}

\newcommand{\ax}{{\rm ax}}
\newcommand{\Acc}{{\rm Acc}}
\newcommand{\Act}{{\rm Act}}
\newcommand{\ded}{{\rm ded}}
\newcommand{\Gm}{{$\Gamma_0$}}
\newcommand{\ID}{{${\rm ID}_1^i(\Oc)$}}
\newcommand{\PA}{{\rm PA}}
\newcommand{\ACA}{{${\rm ACA}^i$}}
\newcommand{\RefP}{{${\rm Ref}^*({\rm PA}(P))$}}
\newcommand{\RefS}{{${\rm Ref}^*({\rm S}(P))$}}
\newcommand{\Rfn}{{\rm Rfn}}
\newcommand{\tar}{{\rm Tarski}}
\newcommand{\UNFA}{{${\mathcal U}({\rm NFA})$}}

\author{Nik Weaver}

\title [Intuitionism and the liar paradox]
       {Intuitionism and the liar paradox}

\address {Department of Mathematics\\
          Washington University in Saint Louis\\
          Saint Louis, MO 63130}

\email {nweaver@math.wustl.edu}

\date{\em April 13, 2010}

\begin{abstract}
The concept of informal mathematical proof considered in intuitionism
is apparently vulnerable to a version of the liar paradox. However,
a careful reevaluation of this concept reveals a subtle error whose
correction blocks the contradiction. This leads to a general resolution
of the classical semantic paradoxes.

This paper is an expanded version of parts of \cite{W3}.
\end{abstract}

\maketitle


\section{Informal proof}

\subsection{Validity of informal proofs.}
The concept of a mathematical proof is important in intuitionism. We
think of proofs as objects, and it is supposed to be decidable
whether a given object is or is not a valid proof.

Here ``proof'' does not necessarily mean ``formal proof''.
It is indeed decidable, in the sense of computability theory, whether
a given finite sequence of formulas constitutes a proof within a
given recursive formal system; but that isn't what is meant.
The issue is not syntactic correctness of formal proofs, but
rather semantic correctness of informal proofs.

This probably sounds rather vague. Why not just
agree that ``valid mathematical proof'' means ``proof that can be
formalized in ZFC (Zermelo-Fraenkel set theory, the generally
accepted axiom system for mathematics)''?

\subsection{Going beyond a given formal system.}
Unfortunately, this proposal leads to a dilemma. Do we really know that
the Zermelo-Fraenkel axioms are valid? If so, then we should be able to
infer that ZFC is consistent, i.e., $0 = 1$ is not a theorem. (More
precisely, we should be able to infer a number
theoretic sentence that arithmetically expresses in a standard way
the consistency of ZFC.) But by G\"odel's  second incompleteness theorem,
this is not provable within ZFC (assuming ZFC is, in fact, consistent).
So if we know that the Zermelo-Fraenkel axioms are valid then we can infer,
informally but correctly, a statement that can't be proven within ZFC.
Adding this statement to ZFC then yields a stronger formal system that we
still know to be valid.

The other horn of the dilemma is that if we don't know whether the
Zermelo-Fraenkel axioms are valid, then there are proofs executable in
ZFC whose validity is in question. So either way, we cannot regard
ZFC as exactly capturing all valid informal reasoning.

The point is not specific to ZFC. Any formal system $S$ that
might be proposed as exactly modelling all valid informal reasoning
would be subject to the same complaint. Either we don't know that the
axioms of $S$ are valid, in which case $S$ manifestly fails to achieve
its goal, or we do know that they are valid, in which case we can go
beyond $S$. We can strengthen it by augmenting it with an assertion of
its consistency.

What all this tells us is that we're dealing with a free-floating
concept of ``valid proof'' that not only isn't thought of as being
tied to any particular formal system, it actually {\it can't} be
exactly modelled by any formal
system. The concept is inherently unformalizable. To be sure, it
can be partially formalized: we can produce formal systems that do
capture some aspects of valid informal reasoning --- maybe even, for
practical purposes, all important aspects. But we can always
go beyond any given partial formalization.

\subsection{Is the proof concept meaningful?}
At this point one could ask why we have to accept that
this informal notion of ``valid mathematical proof'' is even coherent.
The fact that it can't be formalized seems like good evidence
that it simply isn't meaningful. But if we had no overarching notion
of validity, all that would be left would be the syntactic notion of
validity internal to a given formal system. We would then have to conclude
that mathematics is nothing but a meaningless formal game with symbols.

Surely most mathematicians believe, for example, that there really are
infinitely many prime numbers, and that Euclid's proof of this fact is not
merely valid in the trivial sense of being syntactically correct within
some formal system, it is valid in a general semantic sense. For this to
be the case we need to have a general semantic notion of proof validity.

\subsection{Indefinite extensibility.}
So, regardless of how we feel about intuitionism,
unless we're going to be hardcore formalists it would seem that we
have to accept as meaningful the informal notion of semantic validity of
mathematical proofs. The problem then becomes, perhaps, a psychological
one: how to come to terms with the fact that we have to deal with
a concept that can only be partially formalized.

Michael Dummett's ideas about ``indefinitely extensible'' concepts
may be helpful here. According to Dummett (\cite{Dum}, p.\ 441) there
is a whole category
of concepts which have a special kind of productive quality that he
calls indefinite extensibility. By this he means that whenever we've
precisely circumscribed some definite collection of individuals falling
under such a concept, it will always be possible to find a new
individual falling under the concept that wasn't captured.
This is exactly what we just observed about informal proofs: any
partial formalization can be extended. So ``valid proof'' is an
indefinitely extensible concept.

But it isn't the only one. Dummett would also say that the concept
``set'' is indefinitely extensible. And we can point to other
familiar concepts that have the same kind of productive quality
(e.g., ``truth'', ``definition''). So if we buy into the idea of indefinitely
extensible concepts then we may come to feel that our inability to
formalize the notion ``valid mathematical proof'' is not a defect
of that notion, but simply an expression of the fact that it's
indefinitely extensible. To put this another way, we ought to be
comfortable with the concept ``valid proof'', despite the fact
that it's only capable of partial formalization, to the same degree
that we're comfortable with the concept ``set'', which is also
only capable of partial formalization. (The truths of first order
set theory are not recursively enumerable.)

How is it possible to reason about indefinitely extensible concepts?
It seems hard to get one's hands on such a slippery notion. Dummett's
answer to this question is that we need to use intuitionistic logic
when we're working with assertions that quantify over an indefinitely
extensible concept. This point is debatable and I will return to it
later in the paper.

\subsection{The provable liar paradox.}
So far, we've done the following. We introduced the informal notion of
``valid proof'', observed that it cannot be formalized, convinced
ourselves that we nonetheless need this concept if we are to avoid
hardcore formalism, and (perhaps) mollified our concerns about its
unformalizability by placing it in the general setting of indefinitely
extensible concepts.

But there's another problem: the informal notion of proof seems to be
directly implicated in paradoxes. Specifically, if we accept that there
is a meaningful free-floating notion of proof that is not attached to any
particular formal system, then it is hard to see why we can't formulate a
sentence that asserts of itself that it is not provable in this free-floating,
informal sense:
$$L = \mbox{``This sentence is not provable.''}$$
This sentence is paradoxical in the same way as the usual liar sentence.
But recall that Dummett tells us to use intuitionistic logic, in which
the law of excluded middle ($A \vee \neg A$, for all formulas $A$) is not
assumed. In particular, we are not forced to accept the dichotomy ``either
$L$ is provable or $L$ is not provable''; maybe this is the key to evading
the paradox? Unfortunately, no. It is easy
to see that we can deduce a contradiction using pure intuitionistic logic,
without invoking excluded middle. (First check that assuming
$L$ is provable leads to a contradiction, then infer from this that $L$
is not provable. This is a proof that $L$ is not provable, and hence it
is a proof of $L$. So $L$ is provable, which leads to a contradiction.)

To my mind this is a very serious problem for the informal notion of proof.
The direct involvement of this notion in a paradox seems like an excellent
reason to suppose that it is not meaningful.

(True, the informal ``set'' concept is also apparently paradoxical. But
this paradox is resolved by recognizing that the ``set'' concept is
indefinitely
extensible. The general principle would be that we can form the set of
$x$ such that $P(x)$ only for ordinary predicates $P$, not indefinitely
extensible ones. In the specific case of Russell's paradox, for example,
there is no such thing as the set of all sets that aren't members
of themselves, precisely because the concept ``set that is not a member
of itself'' is indefinitely extensible. There is no analogous
resolution of the liar paradox.)

\section{Proof as a heuristic concept}

\subsection{A mistake.}
I'm not going to mount a defense of the intuitionist's proof concept against
this new attack. I agree with the criticism. I think that the informal notion
of valid mathematical proof, as it is characterized in traditional
intuitionism, is indeed incoherent. The proof of this is the fact that
it leads directly to paradox.

What I want to do instead is, first, to identify what I think is a mistake
intuitionists have made, and second, to show how correcting this mistake
defuses the paradox. Then I'll conclude by making a case that this also
settles the usual liar paradox which is based on truth rather than
provability.

The intuitionist's mistake is something I mentioned at the beginning of the
paper, the assumption that proof validity is decidable. This seems to me to
run firmly against general intuitionistic ideas. Consider that intuitionists
would say that we cannot, at present, affirm that Goldbach's conjecture has a
definite truth value. In order to do this, they would argue,
we would need either (1) a proof of
the conjecture, (2) a counterexample, or (3) at a minimum, a procedure which
is guaranteed to terminate after a finite number of steps and produce one of
the first two items. If we don't have a finite procedure that will produce
either a proof or a counterexample, we can't assert that the conjecture
has a definite truth value.

Now suppose we are given what appears to be an informal proof
of Goldbach's conjecture. How would we check it? Do we have a finite
procedure that will, in principle, infallibly assess the validity of
any potential informal proof? If not, then it would seem that we can't
affirm that the statement ``$p$ is a valid proof'' has a definite truth
value for all $p$ without violating intuitionistic principles.

Of course, if we're talking about formal proofs, then we do have such a
procedure. But the whole point of considering informal proofs is that this
takes us beyond formal proofs in any specified formal system. Decidability
becomes an issue if the purported proof requires not just the axioms of
some accepted base system (Peano arithmetic, say) but additional principles
which might or might not informally follow in some way from the accepted
axioms. For instance, along the lines we discussed above, if we recognize
the Peano axioms as valid then we can accept not only Peano arithmetic (\PA),
but also the stronger system
$$\PA' = \PA + {\rm Con}(\PA),$$
where ${\rm Con}(\PA)$ is a number theoretic sentence that arithmetically
expresses in a standard way the consistency of PA. Then we can go one
step further and accept
$$\PA'' = \PA' + {\rm Con}(\PA'),$$
and so on. This process can be iterated any finite number of steps,
and (with minor technical complication) even beyond that into the transfinite
\cite{Tur, Fef}. (In particular, a simple informal induction argument can
lead us to accept $\PA^{(n)}$ for all $n$.) Just how far we can go is a
subtle question, and it apparently cannot be decidable, as this would
contradict indefinite extensibility of the ``proof'' concept. If we could
algorithmically decide which systems $\PA^{(\alpha)}$ we can accept, then
we could formulate a single system whose theorems are precisely the
theorems of all the acceptable $\PA^{(\alpha)}$, and we could then go
beyond it.

Troeltra (\cite{Tro}, p.\ 7) suggests that we can, in effect,
build decidability into our
informal notion of proof, by stipulating that $p$ does not count as a proof
unless we have no doubt that it is a proof. The problem with this idea
is that it assumes we can decide whether there is any doubt about whether
$p$ is a proof. In other words, it simply shifts the assumption of
decidability onto a different predicate, where it is equally unfounded.

\subsection{Heuristic concepts.}
Let's grant that proof validity is not decidable and see what the
consequences are.

First, we need a new name. We can no longer say that the ``valid
proof'' concept
is indefinitely extensible since this term, in Dummett's usage, implies
decidability. I will say that the informal concept of proof is {\it heuristic}
\cite{W1}. The point is that if $C$ is an indefinitely extensible
concept then we can always produce new individuals that fall under $C$
by going beyond the current repertoire of available individuals,
but our understanding of what it means to fall under $C$ does not change.
If $C$ is heuristic, on the other hand, then one way to enlarge the pool
of individuals falling under $C$ is to induct previously available
individuals by expanding our conception of what it means to fall under
$C$. In the same way that Goldbach's conjecture could change its
status from ``not known to have a truth value'' to, say, ``true'', an
existing potential proof of Goldbach's conjecture could change its status
from ``not known to be valid or invalid'' to, say, ``valid''.

We've already balked at the prospect of dealing with concepts that
can't be fully formalized, but were at least decidable. Now we're
going one step further. Is it possible to reason in any way with
heuristic concepts, or is this category just too vague to be dealt
with at all?

\subsection{Natural deduction.}
It's easiest to get at this question using natural deduction. This is a
system of deduction rules which directly express the meanings of the
logical symbols (see, e.g., \cite{Pra}). For instance, one rule states
that given $A$ we can deduce $A \vee B$.

We adopt the intuitionistic ``proof interpretation'' of the logical symbols.
This means that
a proof of $A \wedge B$ is by definition a proof of $A$ together with a
proof of $B$, a proof of $A \vee B$ is either a proof of $A$ or a proof of
$B$, and so on. Under this interpretation, is it fair to deduce $A \vee B$
from $A$? That is, if there is a proof of $A$, can we legitimately infer
that there is a proof of $A \vee B$, i.e., a proof of $A$ or a proof of $B$?

Clearly the answer is: yes, this is a legitimate inference. The
possibility that there may be cases where we cannot decide whether
we have a proof of $A$ is irrelevant to the question since the premise
of this inference is that there definitely is a proof of $A$.

At the risk of belaboring the point, we could make a similar inference
for any concept. If there is a wixle of $A$ then we may correctly infer
that there is a wixle of $A$ or a wixle of $B$. As long as ``wixle'' is
a meaningful concept, this is legitimate regardless of what ``wixle''
actually means, and even if the predicate ``is a wixle'' is not decidable.

What about the other rules of natural deduction? For example, the modus
ponens rule states that we may infer $B$ from the two hypotheses $A$ and
$A \to B$. Under the proof interpretation of the logical symbols,
a proof of $A \to B$ is a procedure that converts any proof of $A$ into
a proof of $B$. So is modus ponens justified if proof is heuristic?  Yes.
In general, if there is a wixle of $A$ and there is a procedure that converts
any wixle of $A$ into a wixle of $B$, then there is a wixle of $B$. Once
again, this in no way makes use of any assumption about decidability
of the ``wixle'' concept.

It begins to look as though dropping the assumption of decidability
does not have any dramatic effects.

\subsection{Circularity.}
There is one important consequence, however: if the ``proof'' concept is
heuristic then we cannot proceed from the assumption that all proofs are valid.

Suppose ``proof'' is decidable. That is, it is a completely
definite concept all of whose properties are determined in advance.
One of these properties is that anything that counts as a proof must
be valid. So as we identify and adopt various axioms and deduction rules,
we can legitimately make use of the fact that all proofs establish correct
conclusions. We may use this fact to help ourselves identify which
axioms and rules are valid, if it is of any help.

If ``proof'' is merely heuristic, on the other hand, then our adoption of
axioms and deduction rules should be seen in a very different light. We
are not merely identifying the properties of an already definite ``proof''
concept, we are building up the concept itself. In this case it is not
reasonable to assume that all proofs will ultimately turn out to be valid.
If that assumption in any way feeds into our choice of which axioms and
rules to adopt, it would be circular because we could end up adopting an
axiom or rule for reasons which hinge on the correctness of proofs that
themselves use the axiom or rule in question.

An analogy may help here. Suppose we accept the validity of Peano arithmetic.
Then we can also accept the stronger system
$$\PA' = \PA + {\rm Con}(\PA).$$
But what about a system like
$$\PA^* = \PA + {\rm Con}(\PA^*)?$$
It is easy, using G\"odellian self-reference techniques, to formulate
a sentence ${\rm Con}(\PA^*)$ that arithmetically expresses the consistency
of PA augmented by the sentence ${\rm Con}(\PA^*)$ itself. However, this
stronger assertion is incorrect. Since $\PA^*$ proves its own
consistency, we know from the second incompleteness theorem that it is
inconsistent. Thus ${\rm Con}(\PA^*)$ can be disproven within PA.

What went wrong? We are free to adopt a new principle that acknowledges
the consistency of previously adopted principles. But in general we cannot
adopt a new principle whose correctness hinges on the correctness of the
system as a whole, {\it including the new principle that is being adopted}.
That would be circular.

Again, if ``proof'' were decidable then this would not be an issue. There
would be nothing wrong in assuming the global correctness of the ``proof''
concept, if doing so would help us to identify which proof principles
are valid. But if ``proof'' is heuristic then there is a very real
danger of circularity in such an assumption. We cannot adopt any axiom
or deduction rule whose justification implicitly assumes the global validity
of all proofs, including proofs that might make use of the new axiom or rule
under consideration.

\subsection{Ex falso quodlibet.}
How could such an assumption arise? We can give a simple example involving
the ``ex falso quodlibet'' law which states that anything follows from a
contradiction.

Before I make this point, let me acknowledge that there are many settings
in which the ex falso law is legitimate. It is commonly accepted, for
example, that in the setting of first order arithmetic we can give a
direct argument that any formula follows from $0=1$.

I have elsewhere argued that this fact generalizes to any system that
proves the law of excluded middle for all atomic formulas; see Section
2.3 of \cite{W2}. That is, ex falso is valid if all of the basic concepts
in play are decidable.

But is it universally valid? There are several possible ways one could try
to justify this conclusion. I want to consider the following argument,
which I will call the {\it justification by vacuity}:
\begin{quote}
``For any formula $A$, the assertion $0=1$ $\to$ $A$ means that there is a
procedure which will convert any proof of $0=1$ into a proof of $A$. But
there are no proofs of $0=1$, so this is vacuously the case: for any
formula $A$, the null procedure will convert any proof of
$0=1$ into a proof of $A$. Thus $0=1$ $\to$ $A$.''
\end{quote}
(This is paraphrased from \cite{vD}.) On its face, the argument is persuasive.
The catch is that in order for it to work we must know that there is no
proof of $0=1$, and this hinges on the global correctness of all proofs.
The justification by vacuity is circular because its justification of the
ex falso law assumes that no proof, {\it including proofs that might make
use of the ex falso law}, establishes $0=1$.

There may be some other way to universally justify ex falso; I don't think
so, but that is not important here. The only
point I want to make is that the justification by
vacuity exhibits the kind of circularity that we have identified as
illegitimate. It is a textbook example of an attempted justification
which is circular in the manner discussed above.

\subsection{An objection.}
Readers of \cite{W3}, where this analysis was first presented,
have objected that the circularity I identify in ex falso is no worse than
related, or possibly identical, circularities which are present in all
of the standard deduction rules. For example, in some sense
the rule ``given $A$, infer $A \vee B$'' is circular
because a proof which establishes $A$ may itself have employed
this rule. So don't we need to assume the global correctness of all proofs
in order to justify this rule too? And couldn't the same be said of any
of the rules of natural deduction?

No. First of all, the alleged circularity present in the rule
``infer $A \vee B$ from $A$'' is not the same as the circularity
we have identified in the justification of ex falso discussed above. This can
easily be seen by applying the ``wixle'' test: would either justification
still be valid if we replaced ``proof'' with ``wixle'', where ``wixle'' could
be any meaningful concept? Consider:
\begin{enumerate}
\item
if there is a wixle of $A$ then there is a wixle of $A$ or a wixle of $B$
\item
any wixle of $0=1$ can be converted into a wixle of $A$.
\end{enumerate}
We don't need to know anything about wixles to be sure that if there
is a wixle of $A$ then there is a wixle of $A$ or a wixle of $B$. But we
do need to know something about wixles to be sure that any wixle of $0=1$
can be converted into a wixle of $A$, for arbitrary $A$. Namely, in order for
the justification by vacuity to work, we would need to know that there
are no wixles of $0=1$. If ``wixle'' means ``proof'', then this means that
we require the global correctness (or at least consistency) of all proofs.
This is where the circularity comes in. There is a circularity in the
justification of (2) which is {\it not present} in the justification of (1).

In order to clarify this distinction, we need to understand better just
what is wrong with circularity. The danger to avoid is a justification of
$X$ that implicitly assumes that $X$ is correct. This would beg the question,
and such a justification is not to be trusted. That is why we must reject any
justification of a proof technique that assumes the global correctness of
all proofs. If the justification succeeds, then ``all proofs'' would include
proofs that make use of the technique in question, so that the assumption
of global correctness of proofs would have presumed the correctness of the
technique which was supposedly being justified. This is how circularity can
be dangerous.

Given an attempt to justify a proof principle $X$, a simple way to test for
this kind of circularity is to ask the question, ``If $X$ were not correct,
but we nonetheless admitted proofs which made use of $X$,
would this affect the validity of the justification?'' If it would, then
the justification is dangerously circular; if not, then the justification
is not circular in any critical way.

So suppose the inference ``$A$ entails $A \vee B$'' were not generally
valid, but we nonetheless admitted proofs that made use of this inference.
Would this affect the correctness of the reasoning that if there is a proof
of $A$ then there is a proof of $A$ or a proof of $B$? No, that trivial
reasoning would still be valid. ``If there is a wixle of $A$ then there is
either a wixle of $A$ or a wixle of $B$'' works just as well if ``wixle''
is a defective ``proof'' concept. Now suppose the inference ``$0=1$ entails
$A$'' were not generally valid, but we nonetheless admitted proofs that
made use of this inference.
Would this affect the correctness of the justification by vacuity?
Yes, it would, because proofs which made use of ex falso might not
be correct, and hence could conceivably establish $0=1$. But the existence
of proofs of $0=1$ would invalidate the justification by vacuity.
I do not see how I can make it any clearer that justification by vacuity
requires the global validity of proofs in a way that the justifications of
the other deduction rules do not.

\subsection{Modus ponens.}
A more interesting objection that I would rather have seen specifically
targets modus ponens. Given a proof of $A$ and a procedure that converts
any proof of $A$ into a proof of $B$, we may obtain a proof of $B$ ---
but only if the procedure works. There could be a circularity issue here.
Perhaps the concept ``procedure that successfully converts proofs into
proofs'' is heuristic. Moreover, perhaps we can make use of valid proofs
to help construct successful proof conversion procedures. Then the validity
of modus ponens would require that all proof conversion procedures succeed,
which would in turn require the validity of all proofs, including proofs
that made use of modus ponens.

This is a plausible objection, but it is not critical because the variety of
proof conversion procedures that are actually needed in typical formal systems
is quite limited. For instance, one such procedure is ``given a proof of $A$
and a proof of $B$ that assumes $A$, append the latter to the former to produce
a proof of $B$''. If we possess a proof of $B$ assuming $A$, we can use this
procedure to convert any proof of $A$ into a proof of $B$. This is what
justifies the natural deduction rule which allows us to deduce $A \to B$ from
the existence of a proof of $B$ that assumes $A$. It could also happen that
the justification of some nonlogical axiom involving implication might
require the use of some other more special kind of proof conversion
procedure. For instance, in order to justify transitivity of equality,
$$x = y\quad \wedge\quad y = z\quad \to\quad x = z,$$
we need a procedure which will convert a proof of $x = y$ together with
a proof of $y = z$ into a proof of $x = z$. One such procedure is: first
give the proof of $x = y$, then give the proof of $y = z$ with every
occurence of $y$ replaced by $x$.

Other nonlogical axioms might require other special proof conversion
techniques. But generally speaking, in order to justify the axioms of any
particular formal system we would never need to invoke an open-ended
notion of ``proof conversion procedure''. Modus ponens could always
be saved by restricting $A \to B$ to mean ``there is a procedure [of some
definite type] which will convert any proof of $A$ into a proof of $B$''.

\subsection{Minimal logic.}
The basic rules of natural deduction for {\it minimal logic} are:
\begin{enumerate}
\item
given $A$ and $B$, deduce $A \wedge B$;
\item
given $A \wedge B$, deduce $A$ and $B$;
\item
given either $A$ or $B$, deduce $A \vee B$;
\item
given $A \vee B$, a proof of $C$ from $A$, and a proof of $C$ from
$B$, deduce $C$;
\item
given a proof of $B$ from $A$, deduce $A \to B$;
\item
given $A$ and $A \to B$, deduce $B$.
\end{enumerate}
(There are also rules for quantifiers; see \cite{Pra} for a fuller and
more precise account.) We have no special axioms for negation because we
interpret $\neg A$ to mean $A \to 0=1$. {\it Intuitionistic logic} is
minimal logic plus ex falso; {\it classical logic} is intuitionistic
logic plus excluded middle.

From what we have said above, minimal logic is suitable for reasoning
about heuristic concepts, in particular for reasoning about provability.
Ex falso might also be justified in such a setting, but it is not clear that
this will always be the case.

Excluded middle is generally not a valid assumption when heuristic concepts
are in play. This conclusion is already argued by Dummett when one has
to deal with assertions which quantify over indefinitely extensible
concepts. But the problem for heuristic concepts is much more severe.
Indeed, in the indefinitely extensible case, the concepts are assumed to
be decidable and so we do have excluded middle for all atomic formulas.
The only question is whether we can accept it for formulas that quantify
over an indefinitely extensible concept. This would appear to depend on how
we interpret truth: under the proof interpretation, excluded middle is dubious
because there is no reason to assume that every such formula is in principle
provable or disprovable. Since the concept is indefinitely extensible there is
no question of being able, even in principle, to verify such formulas
in any mechanical way. The only way to assess the truth of formulas
with quantification is by deductive reasoning, and generally speaking it
seems unlikely that the truth value of every formula could
be determined in this way. So interpreting ``true'' as ``provable in
principle'' probably renders excluded middle invalid.

However, I do not see anything wrong with a classical
interpretation of truth --- taking ``true'' to mean ``is the case'' rather
than ``is provable'' --- in the indefinitely extensible setting, and this
would support the law of excluded middle.
In any case, at a technical level it seems likely that we could generally
justify excluded middle by replacing each indefinitely extensible
concept with a definite subconcept. For instance, in set theory
replace ``set'' by ``set of accessible rank''. Or, if inaccessible
cardinalities are in play, go to some larger cardinal. Broadly speaking,
there should
always be a cutoff that transcends any previously specified method of
construction. Thus, we should generally be able to justify excluded
middle by reinterpreting the system in a definite ``toy'' universe.

Heuristic settings are radically different. Even when working within
a definite ``toy'' universe it may be impossible to pin down exactly
which individuals fall under a given heuristic concept. For example,
consider the heuristic notion of ``definability'' and the purported
definition
\begin{quote}
the smallest natural number not definable in ten English words
\end{quote}
(Berry's paradox). Even restricting to the finite set of all ten word
long strings, we still arrive at a paradox if we suppose that every such
string definitely either does or does not define a natural number.
So assuming the law of excluded middle when reasoning about heuristic
concepts may be not merely unjustified, but actually inconsistent.

\section{The liar paradox}

\subsection{Reasoning about provability.}
We can now resolve the provability version of the liar paradox discussed
earlier. Suppose we want to reason about provability itself. We introduce
a predicate ${\rm Prov}(\ulcorner A\urcorner)$ signifying that there is a
proof of the formula $A$ with G\"odel number $\ulcorner A\urcorner$.
Not a proof within any particular formal system, but a semantically
valid proof in an informal sense.

Since ``proof'' is heuristic, classical logic will not be appropriate here.
Ex falso might be justifiable, but excluded middle presumably is not.

Which nonlogical axioms can we adopt? Is it legitimate to assume
$A \leftrightarrow {\rm Prov}(\ulcorner A\urcorner)$? Recall that we are
using the proof interpretation of the logical symbols, so that this means
that any proof of $A$ can be converted into a proof that there is a proof
of $A$, and vice versa. One implication is valid: given a proof $p$ of $A$,
we can prove that there is a proof of $A$ by exhibiting $p$. If $p$ really
is a proof of $A$ then it should be possible to verify this. In other words,
we reject the possibility of a valid proof whose validity, in principle,
could never be recognized. Thus, we can accept the axiom
$$A \to {\rm Prov}(\ulcorner A\urcorner).$$

In the reverse direction, can we convert any proof that $A$ has a proof
into a proof of $A$? Here is an argument suggesting that we can. Since we
are reasoning constructively, any proof that some object exists ought to,
in principle, actually give us a means of constructing that object. So any
proof that there is a proof of $A$ should, if executed, actually produce
a proof of $A$. This means that we can convert a proof $p$ that $A$ is
provable into a proof of $A$ by executing $p$ and displaying the result.

But this argument is flawed in the same way that the justification by
vacuity of ex falso is flawed. It requires that any proof that $A$ has a
proof, {\it including proofs that might make use of the axiom
${\rm Prov}(\ulcorner A\urcorner) \to A$}, actually does produce a proof of
$A$. That is, it assumes the global correctness of all proofs. So we
{\it cannot} adopt the axiom ${\rm Prov}(\ulcorner A\urcorner) \to A$.

Additional axioms for provability that are straightforwardly justifiable
include
\begin{eqnarray*}
{\rm Prov}(\ulcorner A\urcorner) \wedge {\rm Prov}(\ulcorner B\urcorner)
&\leftrightarrow& {\rm Prov}(\ulcorner A \wedge B\urcorner)\cr
{\rm Prov}(\ulcorner A\urcorner) \vee {\rm Prov}(\ulcorner B\urcorner)
&\to& {\rm Prov}(\ulcorner A \vee B\urcorner)\cr
{\rm Prov}(\ulcorner A \vee B\urcorner) \wedge
{\rm Prov}(\ulcorner A \to C\urcorner) \wedge
{\rm Prov}(\ulcorner B \to C\urcorner)
&\to& {\rm Prov}(\ulcorner C\urcorner)\cr
{\rm Prov}(\ulcorner A\urcorner) \wedge {\rm Prov}(\ulcorner A\to B\urcorner)
&\to& {\rm Prov}(\ulcorner B\urcorner),
\end{eqnarray*}
for all formulas $A$, $B$, and $C$. (See \cite{W3} for additional
explanation.)

\subsection{The provable liar.}
Let $L$ be the ``provable liar sentence'' introduced earlier. One way to
formalize it is as $L = {\rm Prov}(\ulcorner \neg L\urcorner)$.
We can reason informally about $L$ as follows. Assume $L$. Then we have
${\rm Prov}(\ulcorner \neg L\urcorner)$ by definition. But also, since $L \to
{\rm Prov}(\ulcorner L\urcorner)$, we have ${\rm Prov}(\ulcorner L\urcorner)$.
Thus, we have ${\rm Prov}(\ulcorner L \wedge \neg L\urcorner)$ and hence
we have ${\rm Prov}(\ulcorner 0=1\urcorner)$. So we have shown
$$L \to {\rm Prov}(\ulcorner 0=1\urcorner).$$

Or suppose we assume $\neg L$. Then the axiom scheme $A \to
{\rm Prov}(\ulcorner A\urcorner)$ allows us to infer
${\rm Prov}(\ulcorner \neg L\urcorner)$, i.e., we can infer $L$. Thus,
we have $L \wedge \neg L$ and hence we have $0=1$. So we have shown
$\neg L \to 0=1$, i.e., $\neg\neg L$.

We have obtained real, informative conclusions about $L$. It entails
that a contradiction is provable, and its absurdity is absurd. (If the
provable liar sentence is formalized instead as $L =
\neg{\rm Prov}(\ulcorner L\urcorner)$, then we obtain similar results:
we can prove $\neg L$
and $\neg\neg{\rm Prov}(\ulcorner L\urcorner)$.) But lacking
both ${\rm Prov}(\ulcorner A\urcorner) \to A$ and $A \vee \neg A$, we
cannot deduce a contradiction. Both restrictions are important: adding
the former as an axiom would allow the inference $L \to 0=1$,
i.e., we could prove
$\neg L$, which as we just saw leads to a contradiction. Adding
the latter as an axiom would yield ${\rm Prov}(\ulcorner 0=1\urcorner)$,
since both $L$ and $\neg L$ entail this conclusion.

Removing the circularity inherent in adopting axioms whose justification
assumes the global correctness of all proofs defeats the provable liar paradox.
This conclusion is made precise in \cite{W3}, where we present a formal
system HT for reasoning about propositions involving a self-applicative
provability predicate. We prove that neither $0=1$ nor
${\rm Prov}(\ulcorner 0=1\urcorner)$ is a theorem of HT.

\subsection{Revenge.}
Attempts to resolve any version of the liar paradox have to deal
with the so-called ``revenge problem'. Typically a resolution that may
seem to handle the original liar type sentence is defeated by a modified
sentence that is engineered for that purpose.

The resolution offered here does not suffer from any revenge problem
because we do not come to any definite conclusion about the sentence $L$.
We neither affirm it, nor deny it, nor assert that it can be neither affirmed
nor denied. Since we lack excluded middle we are not forced to take any of
these positions. We do assert that $\neg L$ is absurd, but lacking excluded
middle, this alone does not entail $L$.

What we can say is that if $L$ were the case then there would be a proof
of $0=1$. But we cannot infer from this that $L$ is not the case,
precisely because we are not in a position to definitely affirm that there
is no proof of $0=1$. We discussed the illegitimacy of this assumption
at some length; now we see that not being able to prove that we reason
consistently, i.e., not being able to prove
$\neg {\rm Prov}(\ulcorner 0=1\urcorner)$, is fortunate. If we could
prove this then we could combine that result with the formula
$$L \to {\rm Prov}(\ulcorner 0=1\urcorner)$$
derived above to infer $\neg L$, and from this obtain a contradiction.

We know that neither $0=1$ nor ${\rm Prov}(\ulcorner 0=1\urcorner)$ is
a theorem of the formal system HT, which is supposed to model valid informal
reasoning about provability. However, this does not ensure that there is no
proof of $0=1$, since HT is, necessarily, only a partial formalization of
valid informal reasoning.

\subsection{The classical liar paradox.}
Does our resolution of the provable liar paradox also apply to the classical
liar paradox framed in terms of truth?

It does if we identify ``true'' with ``provable'', but this is dubious because
provability is heuristic and our basic intuition about truth tells us that it
is completely sharp. A sentence is true if and only if what it asserts is the
case, and provided we're talking about a sentence that has a definite
meaning, that seems like a completely sharp question.

The qualification that the sentence have a definite meaning is a key
point, because {\it until we define the word ``true''} any sentence that
refers to this concept evidently does not have a definite meaning.
It follows that it would be circular to use
$$\hbox{a sentence is true}\quad\Leftrightarrow
\quad\hbox{what it asserts is the case}\eqno{(*)}$$
as a definition of truth for sentences that themselves contain
the word ``true''.

There are two basic ways of dealing with this circularity. One, due to
Tarski, is to introduce a sequence of words $\{$true${}_n\}$ and apply
($*$) inductively, so that a sentence is true${}_n$ if (1) it does not
contain the word ``true${}_k$'' for any $k \geq n$, and (2) what it
asserts is the case. We get around the problem of circularity by using
a heirarchy of truth predicates, each of which evaluates sentences
that only involve truth predicates at lower levels. The other
possibility is Kripke's suggestion to work with a single truth predicate
but use ($*$) to define its range in stages. So at stage $0$ we define what
it means for sentences that don't contain the word ``true'' to be true, at
stage $1$ we define the truth of sentences that at worst refer to the
truth of stage $0$ sentences, and so on. This leads to a general notion
of ``groundedness'' and a reasonable definition of truth for grounded
sentences. The liar sentence
$$\mbox{``This sentence is not true.''}$$
is an example of an ungrounded assertion, where attempting to evaluate its
truth value leads to an infinite loop, rather than terminating after finitely
many stages as would be the case for a grounded assertion.

Both approaches are interesting and valuable, but both leave something to
be desired as a general explication of truth. Tarski's solution doesn't
do justice to our intuition of truth as a unitary notion. In particular,
we have a clear sense that any sentence that is true${}_n$ for some $n$ is,
in fact, really true. Besides, we want to be able to say things like
$$\hbox{The liar sentence is not true${}_n$ for any $n$.}$$
but according to Tarski {\it this} sentence does not have a truth value.
At any rate, it isn't true${}_n$ for any $n$. Kripke's solution is subject
to a similar objection (as Kripke acknowledges; see pp.\ 714-715 of
\cite{Kri}): we want to say
$$\hbox{The liar sentence is not true, indeed it is not even grounded.}$$
but this sentence itself is not grounded and hence has no truth value in
Kripke's scheme.

The preceding should make it clear that truth in a general setting, which includes
consideration of sentences that themselves refer to truth, does not have
the ``absolute'' quality that it does seem to have in limited settings
where it is applied to a restricted class of assertions, each of which
has a preexisting well-defined meaning. Instead, it has a heuristic
quality. Any partial definition of truth can always be extended further
using ($*$).

The truth theories of Tarski and Kripke are unfaithful to our informal
notion of truth, since both deny truth values to assertions which are
intuitively true.
In contrast, equating truth with provability matches our intuition
perfectly in both cases: under this interpretation, we may affirm that
the liar sentence is neither grounded nor true${}_n$ for any $n$, since
we can easily give informal proofs of these facts.

The objection that provability is heuristic while truth is definite is
no longer persuasive. To the objection that intuitionistic provability
is not adequate to account for mathematical truth, we reply that the
substantive weakness of intuitionistic
mathematics has less to do with the use of intuitionistic logic than it
has to do with the traditional intuitionistic rejection of a completed
infinity. The latter is not relevant to our treatment of the liar paradox;
if we accept the idea of a completed infinity then we may also accept the
possibility of infinite proofs. All that matters is that we have some
coherent notion of provability and that it is heuristic. It need not be
finitary in any sense.

At any rate, if we equate truth with provability then the standard liar
paradox will have the same satisfying resolution as the provable liar
paradox. It is satisfying because it is philosophically well-motivated,
it allows us to draw interesting conclusions about the liar sentence
($L \to {\rm Prov}(\ulcorner 0=1\urcorner)$ and $\neg\neg L$), and it is
technically substantive (HT does not prove $0=1$ or
${\rm Prov}(\ulcorner 0=1\urcorner)$, as we show in \cite{W3}). Moreover,
as explained in \cite{W3}, we can give related resolutions of the paradoxes
of Berry and Grelling-Nelson.

On the other hand, we are free to define ``true'' however we like and there
are certainly settings in which Tarskian or Kripkean definitions are
desirable. If we adopt either of those definitions then the revenge
problem should perhaps simply be ignored. Although we can informally prove
that the liar sentence is neither grounded nor true${}_n$ for any $n$, this
doesn't matter because truth is not equated with provability.

If that conclusion is unsatisfying, this is a testament to the strength
of our intuition that truth at the broadest level really is equatable
with provability.


\end{document}